\documentclass{article}

\usepackage{amsmath,amssymb,amsthm}
\usepackage{tikz}
\usepackage{color}
\usepackage[toc]{appendix}
\usepackage{graphicx}
\usepackage{fancyhdr}
\usepackage{enumitem}
\usepackage{bbm}
\usepackage{parskip}
\usepackage{float}
\usepackage{chngpage}
\usepackage{calc}
\usepackage{bigints}
\usepackage{array}
\usepackage{booktabs}
\usepackage{rotating}
\usepackage{multirow}
\usepackage{adjustbox}
\usepackage{tabularx}
\usepackage{verbatim}
\usepackage{mathtools}
\usepackage{ragged2e}
\usepackage[makeroom]{cancel}
\usepackage{caption}
\usepackage{hyperref}
\usepackage{caption}
\usepackage{subcaption}
\usepackage{appendix}
\usepackage{pgfplots}
\pgfplotsset{compat=1.16}
\usepackage[english]{babel}
\usepackage{hyphenat}
\usepackage[makeindex]{imakeidx}
\usetikzlibrary{datavisualization}
\usetikzlibrary{matrix}
\usetikzlibrary{datavisualization.formats.functions}

\newtheorem{theorem}{Theorem}[section]

\newtheorem{remark}[theorem]{Remark}

\DeclareMathOperator*{\argmax}{arg\,max}

\makeatletter
\def\section{\@startsection {section}{1}{\z@}{3.25ex plus 1ex minus
		.2ex}{1.5ex plus .2ex}{\large\bf}}
\def\subsection{\@startsection{subsection}{2}{\z@}{3.25ex plus 1ex minus
		.2ex}{1.5ex plus .2ex}{\normalsize\bf}}
\@addtoreset{equation}{section} %%
\makeatother
\title{Non trivial optimal sampling rate for estimating a Lipschitz-continuous function in presence of mean-reverting Ornstein-Uhlenbeck noise}

\author{Enrico Bernardi \footnote{Dipartimento di Scienze Statistiche Paolo Fortunati, Università di Bologna. Bologna, Italy.}\, \footnote{\textbf{e-mail}: enrico.bernardi@unibo.it} \and Alberto Lanconelli \footnotemark[1]\, \footnote{\textbf{e-mail}: alberto.lanconelli2@unibo.it} \and Christopher S. A. Lauria \footnotemark[1]\, \footnote{\textbf{e-mail}: christopher.lauria2@unibo.it} \and Berk Tan Perçin \footnotemark[1]\, \footnote{\textbf{e-mail}: berktan.percin2@unibo.it}}

\begin{document}
	
	\maketitle

	\bigskip
	
	\begin{abstract}
	
		We examine a mean-reverting Ornstein-Uhlenbeck process that perturbs an unknown Lipschitz-continuous drift and aim to estimate the drift's value at a predetermined time horizon by sampling the path of the process. Due to the time varying nature of the drift we propose an estimation procedure that involves an online, time-varying optimization scheme implemented using a stochastic gradient ascent algorithm to maximize the log-likelihood of our observations. The objective of the paper is to investigate the optimal sample size/rate for achieving the minimum mean square distance between our estimator and the true value of the drift.
In this setting we uncover a trade-off between the correlation of the observations, which increases with the sample size, and the dynamic nature of the unknown drift, which is weakened by increasing the frequency of observation. The mean square error is shown to be non monotonic in the sample size, attaining a global minimum whose precise description depends on the parameters that govern the model. In the static case, i.e. when the unknown drift is constant, our method outperforms the arithmetic mean of the observations in highly correlated regimes, despite the latter being a natural candidate estimator. We then compare our online estimator with the global maximum likelihood estimator.        	   
	\end{abstract}
	
	Key words and phrases: Ornstein-Uhlenbeck process, positive correlation, log-likelihood function, stochastic gradient ascent. \\
	
	AMS 2000 classification: 60H10, 62F10, 65C20.\\
	
	\allowdisplaybreaks

\section{Introduction and statement of the main result}

High frequency samples of a noisy continuous time signal may produce strongly correlated observations which have a negative impact on the performance of statistical estimation methods. This phenomenon has been observed in a number of different contexts: \cite{Codling2005} and \cite{Rosser2013} explore the complex nature of sampling effects in tracking the movement of individual cells or animals, \cite{Astrom1969} proposes a statistical counterpart of the sampling rate's lower bound from ordinary spectral analysis of time series related to aliasing, \cite{Ait-Sahalia2005} shows that in presence of market microstructure noise the optimal sampling rate for estimating the variance of log returns through their sum of squares is finite, \cite{Zhu2006} proves that the spacing between observations affects the estimation of the variance of a Fraction Brownian Motion.

The aim of the present paper is to analyse this phenomenon in a particular setting. Our framework is very specific but nevertheless pretty common in several applied settings (see for instance \cite{Benth2008}). \\
Specifically, let $\{X_t\}_{t\in [0,1]}$ be our \emph{observation process} which is assumed to be of the form
\begin{align}\label{model}
X_t:=m^{\star}(t)+Y_t,\quad t\in [0,1],
\end{align}
where 
\begin{itemize}
\item $m^{\star}:[0,1]\to\mathbb{R}$ is an \emph{unknown} Lipschitz-continuous function, i.e. there exists $K>0$ such that
\begin{align}\label{11}
|m^{\star}(t)-m^{\star}(s)|\leq K|t-s|\quad\mbox{ for all }t,s\in [0,1];
\end{align}
\item $\{Y_t\}_{t\geq 0}$ is a mean-reverting Ornstein-Uhlenbeck process, i.e.
\begin{align}
\begin{cases}\label{OU}
dY_t=-\theta Y_tdt+\sigma dB_t,\quad t>0;\\ Y_0\thicksim\mathtt{N}\left(0,\frac{\sigma^2}{2\theta}\right).
\end{cases}
\end{align}
Here, the initial condition $Y_0$ is independent of the Brownian Motion $\{B_t\}_{t\geq 0}$ while $\theta$ and $\sigma$ are positive real parameters. The parameter $\theta$ models the \emph{strength of mean-reversion} towards zero and  $\sigma$ the \emph{amplitude of fluctuation} of the process.
\end{itemize}
Our problem is to estimate $m^{\star}(1)$ with a finite sample of observations of the process $\{X_t\}_{t\in [0,1]}$, sampled at equally spaced instants of time. \\ Let $N\in\mathbb{N}$ and assume we observe the process $\{X_t\}_{t\geq 0}$ in the time interval $[0,1]$ at the discrete times $\left\{\frac{1}{N},\frac{2}{N},...,1\right\}$. Then, according to our model specification we have
\begin{align}\label{12}
X_{\frac{t}{N}}\thicksim\mathtt{N}\left(m_t^{\star},\frac{\sigma^2}{2\theta}\right)\quad\mbox{ and }\quad\mathtt{cov}\left(X_{\frac{t}{N}},X_{\frac{s}{N}}\right)=\frac{\sigma^2}{2\theta}e^{-\frac{\theta}{N}|t-s|}\quad\mbox{ for all $t,s\in\{1,...,N\}$},
\end{align}
where for notational convenience we set
\begin{align*}
m^{\star}_t:=m^{\star}\left(\frac{t}{N}\right),\quad t\in\{1,...,N\}.
\end{align*}
We remark that each of the observations $X_{\frac{1}{N}},X_{\frac{2}{N}}...,X_{1}$ has a different unknown mean by virtue of the dynamic nature of the drift $t\mapsto m^{\star}(t)$, nevertheless, assumption \eqref{11} entails that
\begin{align}\label{L_N}
|m^{\star}_t-m^{\star}_s|\leq \frac{K}{N}|t-s|,\quad\mbox{ for all $t,s\in\{1,...,N\}$}.
\end{align}
Thus large values of $N$ dampen the time-varying nature of the sequence $t\mapsto m^{\star}_{t}$. However, large values of $N$ increase -according to equation \eqref{12}- the correlation between observations. Such contrasting effects caused by increases of $N$ are an important element of the setting, we will showcase the existence of a trade-off between them.

We stress that we are facing a non classical statistical problem: our aim is to estimate $m^{\star}(1)$, or in the previously introduced notation $m^{\star}_N$, which is the mean of $X_1$, and to do so we can only utilize (in addition to the value of $X_1$) the observations $X_{\frac{k}{N}}$, $k\in\{1,...,N-1\}$, which are realisations of random variables with means $m^{\star}_1,...,m^{\star}_{N-1}$, different from $m^{\star}_N$. The Lipschitz-continuity assumption that entails \eqref{L_N} provides information on the difference between unknown means and hence will be crucial in deriving an estimator for $m^{\star}(1)$.\\
To estimate $m^{\star}_N$ we set 
\begin{align}\label{13}
\hat{m}_1:=X_{\frac{1}{N}}\quad\mbox{ and }\quad \hat{m}_t:=(1-\alpha)\hat{m}_{t-1}+\alpha X_{\frac{t}{N}},\quad t\in\{2,...,N\},
\end{align}
where the parameter $\alpha\in [0,1]$ is called \emph{learning rate}. Our suggested estimator for $m^{\star}_N$ is $\hat{m}_N$, defined as solution of the recursive equation \eqref{13}. Notice that $\hat{m}_N$ can be written explicitly as
\begin{align}\label{14}
\hat{m}_N=\beta^{N-1}X_{\frac{1}{N}}+\alpha \sum_{j=2}^{N}\beta^{N-j}X_{\frac{j}{N}}.
\end{align}
Scheme \eqref{13} is known in the econometric literature as \emph{simple exponential smoothing}: it was originally defined in \cite{brown1956exponential}, and then further developed in \cite{brown1959statistical} and \cite{brownfore_pre}. Simple exponential smoothing was originally employed as a way to generate forecasts for time series data by giving more weight to recent observations while exponentially diminishing the influence of older data, but it has also been used to estimate the trend of a time series, see for instance \cite{brockwell2009time}.\\
In our setting it is possible to individuate a statistical justification for the estimation rule \eqref{13}. In fact, if we rewrite equation \eqref{13} for $t\in\{2,...,N\}$ as
\begin{align*}
\hat{m}_t=&\hat{m}_{t-1}+\alpha (X_{\frac{t}{N}}-\hat{m}_{t-1})\nonumber\\
=&\hat{m}_{t-1}+\tilde{\alpha}\frac{d}{d\hat{m}}\ln\left(\frac{1}{\sqrt{\pi\sigma^2/\theta}}e^{-\frac{\left(X_{\frac{t}{N}}-\hat{m}\right)^2}{\sigma^2/\theta}}\right)\Bigg|_{\hat{m}=\hat{m}_{t-1}}
\end{align*}
with $\tilde{\alpha}:=\frac{2\alpha\theta}{\sigma^2}$, one can identify \eqref{13} with a stochastic gradient ascent algorithm where the function to be optimized is the log-likelihood of the $t$-th observation, embedding $\hat{m}_N$ in the framework of maximum likelihood estimation. There is however a crucial difference from the usual optimization framework: the likelihood function one wishes to optimize changes at each time $t$, since the observations have different unknown means, and hence the maximum of those functions is time-varying and not static. \\
Tracking a sequence of optima through time is a problem that has already been considered in a variety of settings, see for instance \cite{Polyak, Popkov2005, Zinkevich, Bubeck,Cao}. Using the nomenclature of \cite{simonetto2020time}, these types of problems have been called  \emph{time-varying} optimization problems since their aim is to find the optimum of an objective function that varies with time. \\
The aforementioned lift of the maximum likelihood principle to a dynamic context has been explored, with varying degrees of generality, in recent papers \cite{Lanconelli2023}, \cite{BLL2023} and \cite{BLL2024}, which show how the scheme defined by equation \eqref{13} (or a suitable variation) is able to track the evolution of a time-varying parameter. \\
We are now ready to state our main result whose proof is given in Section \ref{proof main theorem}.
\begin{theorem}\label{main theorem}
	Let $\hat{m}_N$ be defined by the recursive equation \eqref{13}. Then, $\mathbb{E}[|\hat{m}_1-m^{\star}_1|^2]=\frac{\sigma^2}{2\theta}$ and for $N\geq 2$ we have
	\begin{align}\label{upper bound}
	\mathbb{E}[|\hat{m}_N-m^{\star}_N|^2]\leq&\beta^{2(N-1)}\frac{\sigma^2}{2\theta}+\left(\beta^2\frac{1+\beta}{1-\beta}\frac{K^2}{N^2}+\alpha^2\frac{\sigma^2}{2\theta}\frac{e^{\frac{\theta}{N}}+\beta}{e^{\frac{\theta}{N}}-\beta}\right)\frac{1-\beta^{2(N-1)}}{1-\beta^2}\nonumber\\
	&-2\frac{K^2}{N^2}\frac{\beta^{2N}-\beta^{N+1}}{\alpha^2}+\alpha\frac{\sigma^2}{\theta}\frac{e^{\frac{\theta}{N}}-1}{1-\beta e^{-\frac{\theta}{N}}}\frac{\beta^{2N-1}e^{-\frac{\theta}{N}}-\beta^N e^{-\theta}}{\beta e^{\frac{\theta}{N}}-1};
	\end{align}
	here, to ease the notation, we have set $\beta:=1-\alpha$. Moreover, the previous inequality is sharp: if the unknown drift $t\mapsto m^{\star}(t)$ is linear, then it becomes an equality.	
\end{theorem}
\begin{remark}\quad
	\begin{itemize}
		\item The last ratio in \eqref{upper bound} is well defined when $\beta\neq e^{-\frac{\theta}{N}}$. This constraint stems from the requirement of the geometric sum formula, which necessitates this condition for its application. Nevertheless, the result remains valid without imposing this constraint, provided that one eschews the geometric sum formula and retains the complete summation term.
		\item A simple verification shows that the upper bound \eqref{upper bound} converges, as $N$ tends to infinity, to $\frac{\sigma^2}{2\theta}$ which coincides with $\mathbb{E}[|\hat{m}_1-m^{\star}_1|^2]$. This fact yields a non monotonic behavior of the function $N\mapsto \mathbb{E}[|\hat{m}_N-m^{\star}_N|^2]$.
		\item If in \eqref{upper bound} we Taylor-expand $e^{\frac{\theta}{N}}$ in a neighborhood of zero, and neglect terms of order $O\left(\frac{1}{N^2}\right)$, we obtain  
		\begin{align*}
		\mathbb{E}[|\hat{m}_N-m^{\star}_N|^2]\leq&\frac{1+\frac{\theta}{N(1+\beta)}}{1+\frac{\theta}{N(1-\beta)}}\frac{\sigma^2}{2\theta}.
		\end{align*}
		Since the ratio in front of $\frac{\sigma^2}{2\theta}$ is smaller than one, and converges to one as $N$ tends to infinity, we see that the upper bound \eqref{upper bound} converges, as $N$ tends to infinity, to $\frac{\sigma^2}{2\theta}$ from below. 
Thus the non monotonic function  $N\mapsto \mathbb{E}[|\hat{m}_N-m^{\star}_N|^2]$ possesses a minimum.
	\end{itemize}
\end{remark}	

The rest of the paper is organized as follows: in Section 2 we propose several plots of the upper bound \eqref{upper bound} as a function of $N$ with the aim of illustrating the role of the various parameters of the model on the performance of our estimator.
In Section 3, we focus on the case of a constant drift $m^{\star}$ and compare our method with two classical estimators suitable for investigating static problems: the arithmetic mean and the global maximum likelihood estimator. In the final section, we present the proof of our main theorem along with the derivation of two formulas used in the analysis of the static case.

\section{Numerical illustration}

In this section we visually examine the impact of the various parameters in our model \eqref{model}-\eqref{11}-\eqref{OU} on the upper bound \eqref{upper bound}. In particular, we aim to shed light on the relationship between $K$ and the learning rate $\alpha$, as well as the role of $\theta$, in changing the graph of \eqref{upper bound} as a function of $N$.

\begin{figure}[H]
	\centering
	\begin{subfigure}[c]{0.49\textwidth}
		\includegraphics[width=\linewidth]{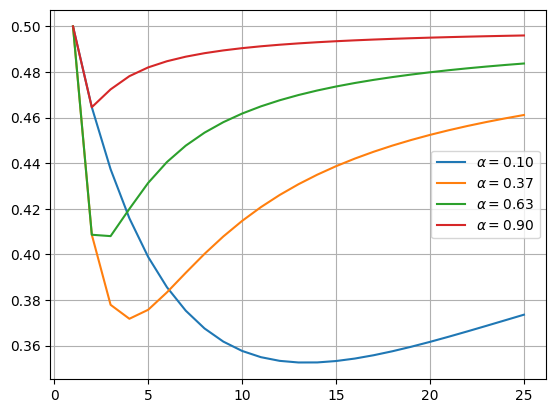}
	\end{subfigure}
	\begin{subfigure}[c]{0.49\textwidth}
		\includegraphics[width=\linewidth]{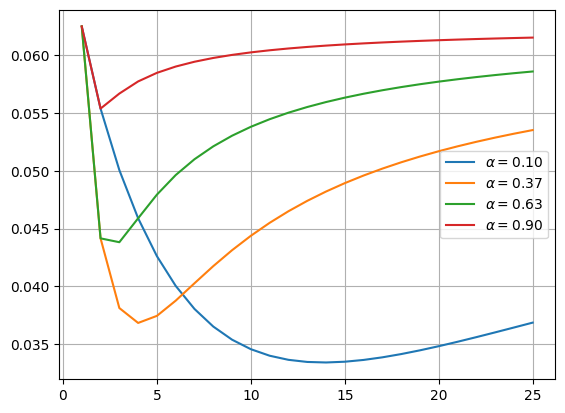}
	\end{subfigure}
\caption{Plots of the sharp upper bound \eqref{upper bound} as a function of $N$ with $K=0$ and different learning rates $\alpha$. In the left figure we set $\theta=1$ and $\sigma=1$ while in the right figure we chose $\theta=2$ and $\sigma=0.5$. The least mean square error is attained in correspondence of small values of the learning rate, as expected in the case of a static drift $m^{\star}$ (i.e. when $K=0$), and for relatively small sample size.} \label{figure 1}
\end{figure}

In Figure \ref{figure 1} we plotted the upper bound \eqref{upper bound} versus $N$ in the static case, i.e. when $K=0$, for two different sets of values for $\theta$ and $\sigma$. We explored four different learning rates, and an inspection of the graphs indicates better performance for lower values of $\alpha$.  

\begin{figure}[h!]
\centering
	\begin{subfigure}[c]{0.49\textwidth}
		\includegraphics[width=\linewidth]{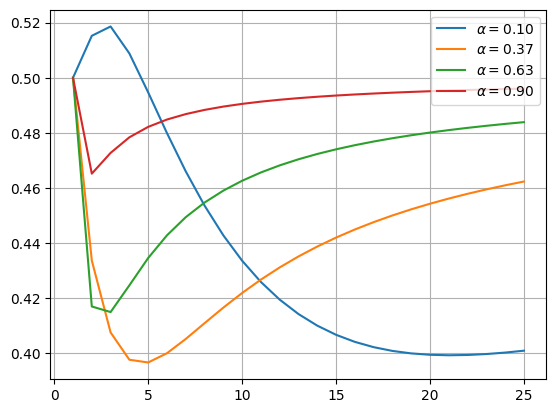}
	\end{subfigure}
	\begin{subfigure}[c]{0.49\textwidth}
		\includegraphics[width=\linewidth]{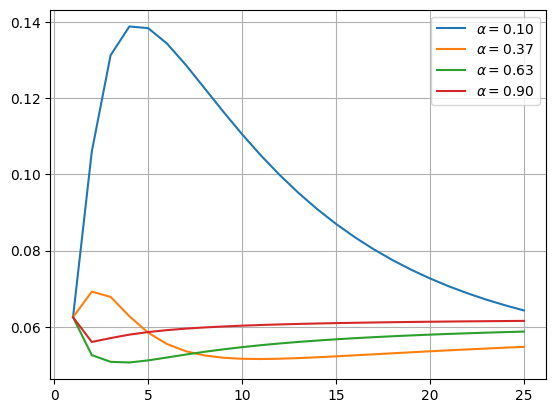}
	\end{subfigure}
\caption{Plots of the sharp upper bound \eqref{upper bound} as a function of $N$ with $K=0.5$ and different learning rates $\alpha$. In the left figure we set $\theta=1$ and $\sigma=1$ while in the right figure we chose $\theta=2$ and $\sigma=0.5$. The increase of $K$ now results in large mean square errors for small learning rates.} \label{figure 3}
\end{figure}

In Figure \ref{figure 3} the Lipschitz constant $K$ is set equal to $0.5$, thus the dynamic nature of the drift $m^{\star}$ assumes importance. In this case the smallest value of the learning rate $\alpha$ does not always result in the best performance.

\begin{figure}[h!]
\centering
	\begin{subfigure}[c]{0.49\textwidth}
		\includegraphics[width=\linewidth]{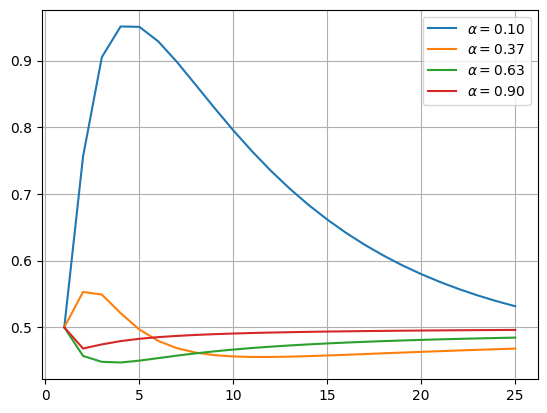}
	\end{subfigure}
	\begin{subfigure}[c]{0.49\textwidth}
		\includegraphics[width=\linewidth]{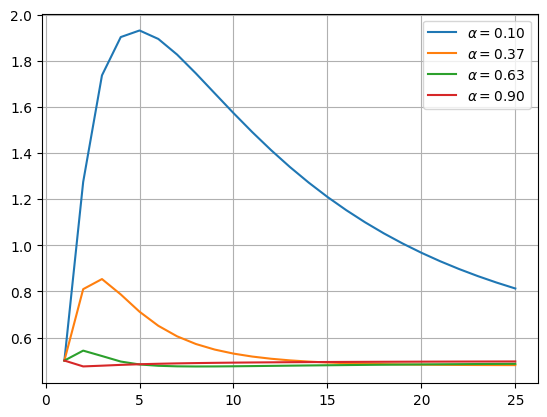}
	\end{subfigure}
\caption{Plots of the sharp upper bound \eqref{upper bound} as a function of $N$ with $K=1.2$ (left figure), $K=2$ (right figure), $\theta=1$, $\sigma=1$  and different learning rates $\alpha$.} \label{figure 5}
\end{figure}

The plots of Figure \ref{figure 5} consider $K=1.2$ (left) and $K=2$ (right). Here, the deterministic behaviour of the drift dominates the uncertainty associated with the observations: the red line in the last plot shows the almost irrelevant contribution of frequent observations in reducing the mean square error. Further increase of $K$ imply large mean square errors in correspondence to small learning rates.

\section{The static case $K=0$} 
		
The case $K=0$, i.e. $m^{\star}(t)=m^{\star}$ for all $t\in [0,1]$, gives us observations $X_\frac{1}{N}$, $X_\frac{2}{N}$,..., $X_{1}$ that are Gaussian with common unknown mean $m^{\star}$ and common variance $\frac{\sigma^2}{2\theta}$. One may therefore wonder whether an application of standard statistical tools can provide lower mean square errors compared to those derived with scheme \eqref{13}. In this section we compare our method with the arithmetic mean of the observations and the (global) maximum likelihood estimator. Here, we emphasize the term 'global' to highlight that our scheme \eqref{13} is also based on the maximum likelihood principle, embedded in the stochastic gradient ascent algorithm, but employed in an online fashion -one observation at a time.
While this feature is unnecessary in the current static framework, it becomes crucial in the broader context presented earlier, where the Lipschitz condition \eqref{11} (and thus the relationship between the unknown parameters $m^{\star}_t$ for $t\in\{1,...,N\}$) implies that the order of observations is significant.

\subsection{Comparison with the arithmetic mean}

\begin{figure}[h!]
	\centering
	\begin{subfigure}[c]{0.49\textwidth}
		\includegraphics[width=\linewidth]{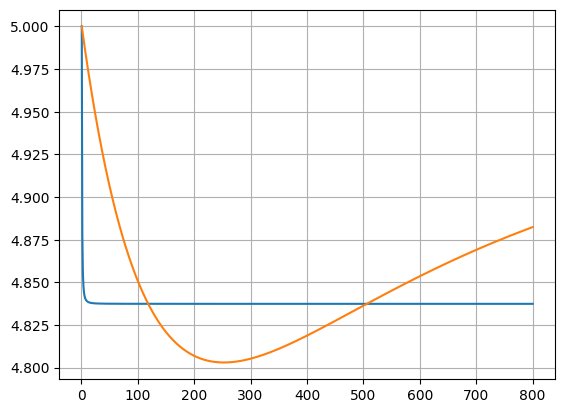}
	\end{subfigure}
	\begin{subfigure}[c]{0.49\textwidth}
		\includegraphics[width=\linewidth]{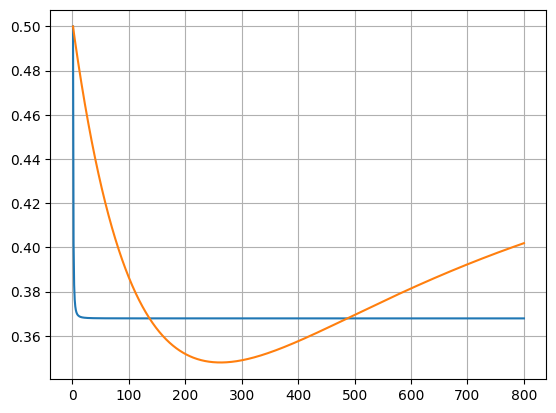}
	\end{subfigure}
\begin{subfigure}[c]{0.49\textwidth}
	\includegraphics[width=\linewidth]{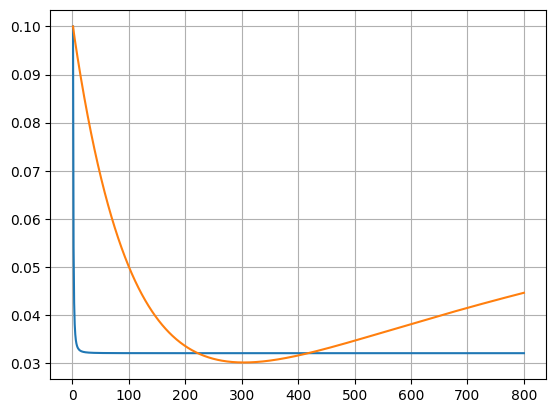}
\end{subfigure}
	\begin{subfigure}[c]{0.49\textwidth}
		\includegraphics[width=\linewidth]{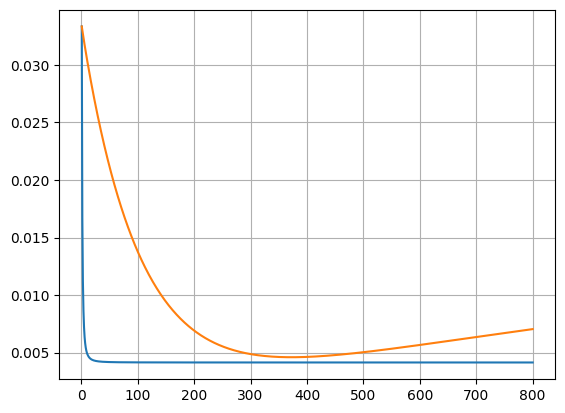}
	\end{subfigure}
\caption{Comparison between the sharp upper bound \eqref{upper bound} with $K=0$ and $\alpha=0.005$ (orange line) and the right hand side of \eqref{LLN} (blue line). The parameter $\sigma$ is set equal to one while $\theta=0.1$ (top left figure), $\theta=1$ (top right figure), $\theta=5$ (bottom left figure) and $\theta=15$ (bottom right figure). The mean square error of our estimator performs better than the arithmetic mean for a suitable range of $N$ values. As $\theta$ increases, that is, when the correlation decreases, the performance disparity between the two statistical procedures diminishes, ultimately resulting in the arithmetic mean surpassing our method. We observe the monotonic behavior of the blue line, in contrast to the orange line, which attains a global minimum at non-trivial values of $N$.}\label{figure 6}
\end{figure}

The mean square error of the arithmetic mean of our observations $X_\frac{1}{N}$, $X_\frac{2}{N}$,..., $X_{1}$ is given by
\begin{align}\label{LLN}
\mathbb{E}\left[\left(\frac{1}{N}\sum_{i=1}^NX_{\frac{i}{N}}-m^{\star}\right)^2\right]=\frac{1}{N^2}\frac{\sigma^2}{2\theta}\frac{Ne^{-\frac{\theta}{N}}-Ne^{\frac{\theta}{N}}+2-2e^{-\theta}}{(1-e^{\frac{\theta}{N}})(1-e^{-\frac{\theta}{N}})},\quad N\geq 1.
\end{align}
A proof of this identity can be found in Section \ref{proof LLN} below. We remark in passing the marginal role of $\sigma^2$ in \eqref{LLN} where it trivially acts as a multiplicative constant.\\
In Figure \ref{figure 6} we compare the performance of \eqref{upper bound} with the one of \eqref{LLN}; the plots describe their behaviour as functions of $N$. It is interesting to note that, if one uses the time-varying learning rate $\frac{1}{t}$ in \eqref{13}, the solution to the corresponding recursive equation coincides with the arithmetic mean  $\frac{1}{N}\sum_{i=1}^NX_{\frac{i}{N}}$. Hence, we are here comparing two different instances of the same general method, one with constant learning rate $\alpha$ and one with vanishing learning rate $\frac{1}{t}$. We remark that, when the drift $m^{\star}$ is non constant, a vanishing learning rate proves inadequate for effectively tracking the dynamic nature of the problem. \\
In Figure \ref{figure 6} we set $\alpha=0.005$, $\sigma=1$ and analyse the sensitivity of the graph to an increase of $\theta$ from $0.1$ up to $15$, i.e., a decrease of the correlation between the observations. When correlation is high, our methods performs better than the arithmetic mean for a specific range of $N$ values; low correlation regimes are instead slightly favorable to the arithmetic mean and provide an almost vanishing mean square error.   

\subsection{Comparison with the maximum likelihood estimator}

\begin{figure}[h!]
	\centering
	\begin{subfigure}[c]{0.49\textwidth}
		\includegraphics[width=\linewidth]{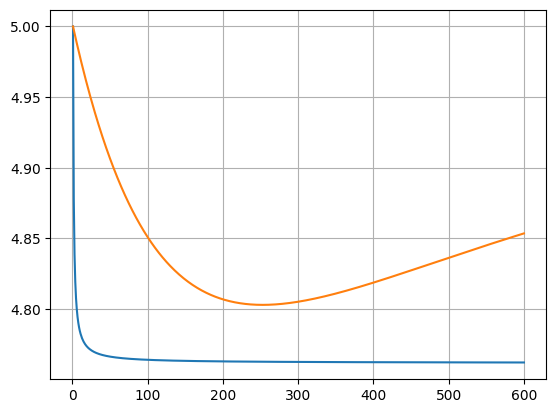}
	\end{subfigure}
	\begin{subfigure}[c]{0.49\textwidth}
		\includegraphics[width=\linewidth]{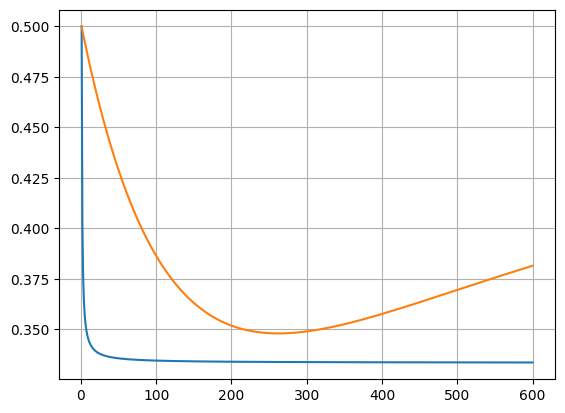}
	\end{subfigure}
	\begin{subfigure}[c]{0.49\textwidth}
		\includegraphics[width=\linewidth]{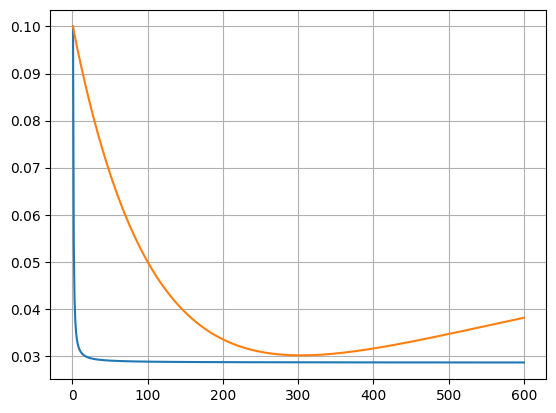}
	\end{subfigure}
	\begin{subfigure}[c]{0.49\textwidth}
		\includegraphics[width=\linewidth]{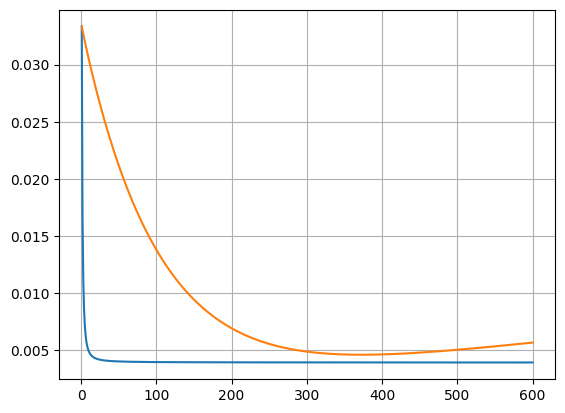}
	\end{subfigure}
	\caption{A comparison between the sharp upper bound \eqref{upper bound} with $K=0$ and $\alpha=0.005$ (orange line) and the right hand side of \eqref{LLN} (blue line). The parameter $\sigma$ is set to one while $\theta=0.1$ (top left figure), $\theta=1$ (top right figure), $\theta=5$ (bottom left figure) and $\theta=15$ (bottom right figure). The (global) maximum likelihood estimator performs better than our method. However, as $\theta$ increases, i.e., as the correlation becomes weaker, the gap between the performances of the two statistical procedures reduces. We again notice the monotonic behavior of the blue line in contrast to the orange one which achieves a global minimum for non trivial values of $N$.}\label{figure 7}
\end{figure}

The (global) maximum likelihood estimator $\hat{M}$ of $m^{\star}$ is defined as
\begin{align*}
\hat{M}:=\argmax_{m\in\mathbb{R}}\frac{1}{(2\pi |C|)^{\frac{N}{2}}}\exp\left\{-\frac{1}{2}\sum_{i,j=1}^Na_{ij}\left(X_{\frac{i}{N}}-m\right)\left(X_{\frac{j}{N}}-m\right)\right\}
\end{align*}
where $C$ stands for the covariance matrix of the vector $\left(X_{\frac{1}{N}},X_{\frac{2}{N}},...,X_{1}\right)$, i.e., $c_{ij}=\frac{\sigma^2}{2\theta}e^{-\frac{\theta}{N}|i-j|}$ for $i,j\in\{1,...,N\}$, $|C|$ denotes the determinant of $C$ while $A=\{a_{ij}\}_{i,j=1,..,N}$ stands for the inverse of $C$.
To compute the $\argmax$ we take the logarithm of the likelihood function and differentiate with respect to $m$, this yields
\begin{align*}
\hat{M}=\frac{\sum\limits_{i,j=1}^Na_{ij}X_{\frac{j}{N}}}{\sum\limits_{i,j=1}^Na_{ij}},
\end{align*}
and its mean square error is given by 
\begin{align}\label{MLE}
\mathbb{E}\left[\left(\hat{M}-m^{\star}\right)^2\right]=\frac{\sigma^2(1-e^{-2\frac{\theta}{N}})}{2\theta\left(N+(N-2)e^{-2\frac{\theta}{N}}-2(N-1)e^{-\frac{\theta}{N}}\right)}.
\end{align}
A proof of this formula can be found in Section \ref{proof MLE}.\\
In Figure \ref{figure 7} we compare the performance of the upper bound of our estimator \eqref{upper bound} with the maximum likelihood one \eqref{MLE}, the plots describe their behaviour as functions of $N$.
Both methods employ the maximum likelihood principle; however, the method associated with scheme \eqref{13} processes one observation at a time, whereas the global maximum likelihood estimator utilizes the entire dataset simultaneously. Consequently, we cannot expect our estimator $\hat{m}_N$ to outperform $\hat{M}$ in the static setting. Nonetheless, the subsequent plots will demonstrate that in low-correlation regimes, the two methods yield nearly equivalent results.
As was the case for the arithmetic mean, in the dynamic setting with the Lipschitz Assumption \ref{11}, the global maximum likelihood estimator loses its efficacy and adapting it to this new context is not straightforward. \\
In Figure \ref{figure 7} we set $\alpha=0.005$, $\sigma=1$ and analyse the sensitivity of the maximum likelihood estimator and our proposed estimator to an increase of $\theta$ from $0.1$ up to $15$, i.e., a decrease of the correlation between the observations.   

\clearpage

\section{Proofs} 

This section is devoted to the proofs of Theorem \eqref{main theorem}, formula \eqref{LLN} and formula \eqref{MLE}.

\subsection{Proof of Theorem \ref{main theorem}}\label{proof main theorem}

The proof given here is inspired by the one contained in \cite{BLL2024}. First of all we recall that 
\begin{align}\label{3}
\hat{m}_t=\beta^{t-1}X_{\frac{1}{N}}+\alpha \sum_{j=2}^{t}\beta^{t-j}X_{\frac{j}{N}}.
\end{align}
Now, 
\begin{align}\label{a}
\hat{m}_t-m^{\star}_t=\beta(\hat{m}_{t-1}-m^{\star}_t)+\alpha \left(X_{\frac{t}{N}}-m^{\star}_t\right),
\end{align}
and taking squares we obtain
\begin{align}\label{2}
(\hat{m}_t-m^{\star}_t)^2=\beta^2(\hat{m}_{t-1}-m^{\star}_t)^2+\alpha^2 \left(X_{\frac{t}{N}}-m^{\star}_t\right)^2+2\alpha\beta(\hat{m}_{t-1}-m^{\star}_t)\left(X_{\frac{t}{N}}-m^{\star}_t\right).
\end{align}
Considering the last term we see that
\begin{align*}
(\hat{m}_{t-1}-m^{\star}_t)\left(X_{\frac{t}{N}}-m^{\star}_t\right)=&\left(\beta^{t-2}X_{\frac{1}{N}}+\alpha \sum_{j=2}^{t-1}\beta^{t-1-j}X_{\frac{j}{N}}-m^{\star}_t\right)\left(X_{\frac{t}{N}}-m^{\star}_t\right)\\
=&\beta^{t-2}X_{\frac{1}{N}}\left(X_{\frac{t}{N}}-m^{\star}_t\right)+\alpha \sum_{j=2}^{t-1}\beta^{t-1-j}X_{\frac{j}{N}}\left(X_{\frac{t}{N}}-m^{\star}_t\right)\\
&-m^{\star}_t\left(X_{\frac{t}{N}}-m^{\star}_t\right)\\
=&\beta^{t-2}\left(X_{\frac{1}{N}}-m^{\star}_1\right)\left(X_{\frac{t}{N}}-m^{\star}_t\right)\\
&+\alpha \sum_{j=2}^{t-1}\beta^{t-1-j}\left(X_{\frac{j}{N}}-m^{\star}_{j}\right)\left(X_{\frac{t}{N}}-m^{\star}_t\right)\\
&+\quad\mbox{zero mean terms}.
\end{align*}
Therefore, taking the expectation in \eqref{2} gives
\begin{align*}
\mathbb{E}[(\hat{m}_t-m^{\star}_t)^2]=&\beta^2\mathbb{E}[(\hat{m}_{t-1}-m^{\star}_t)^2]+\alpha^2\frac{\sigma^2}{2\theta}+2\alpha\beta\left(\beta^{t-2}\gamma(t-1)+\alpha\sum_{j=2}^{t-1}\beta^{t-1-j}\gamma(t-j)\right)\\
=&\beta^2\mathbb{E}[(\hat{m}_{t-1}-m^{\star}_t)^2]+\alpha^2\frac{\sigma^2}{2\theta}+2\alpha\beta\left(\beta^{t-1}\gamma(t-1)+\alpha\sum_{j=1}^{t-1}\beta^{t-1-j}\gamma(t-j)\right)\\
=&\beta^2\mathbb{E}[(\hat{m}_{t-1}-m^{\star}_t)^2]+\alpha^2\frac{\sigma^2}{2\theta}+2\alpha\beta^{t}\gamma(t-1)+2\alpha^2\sum_{j=1}^{t-1}\beta^{t-j}\gamma(t-j)\\
=&\beta^2\mathbb{E}[(\hat{m}_{t-1}-m^{\star}_t)^2]+\alpha^2\frac{\sigma^2}{2\theta}+2\alpha\beta^{t}\gamma(t-1)+2\alpha^2\sum_{k=1}^{t-1}\beta^{k}\gamma(k).
\end{align*} 
Now,
\begin{align*}
\mathbb{E}[(\hat{m}_{t-1}-m^{\star}_t)^2]=&\mathbb{E}[(\hat{m}_{t-1}-m^{\star}_{t-1})^2]+K_t^2-2K_t\mathbb{E}[\hat{m}_{t-1}-m^{\star}_{t-1}],
\end{align*}
with $K_t:=m_t^{\star}-m_{t-1}^{\star}$. On the other hand, taking the expectation in \eqref{a} we get
\begin{align*}
\mathbb{E}[\hat{m}_t-m^{\star}_t]=\beta\mathbb{E}[\hat{m}_{t-1}-m^{\star}_{t-1}]-\beta K_t,
\end{align*}
with $\mathbb{E}[\hat{m}_1-m^{\star}_1]=\mathbb{E}\left[X_{\frac{1}{N}}-m^{\star}_1\right]=0$; this means
\begin{align*}
\mathbb{E}[\hat{m}_t-m^{\star}_t]=-\beta\sum_{j=2}^{t}\beta^{t-j}K_j.
\end{align*}
Therefore, setting $D_t:=\mathbb{E}[(\hat{m}_t-m^{\star}_t)^2]$ we obtain
\begin{align}\label{1}
D_t=&\beta^2\left(D_{t-1}+K_t^2+2 K_t\sum_{j=2}^{t-1}\beta^{t-j}K_j\right)+\alpha^2\frac{\sigma^2}{2\theta}+2\alpha\beta^{t}\gamma(t-1)+2\alpha^2\sum_{k=1}^{t-1}\beta^{k}\gamma(k).
\end{align}
At this point we upper bound each $K_j$ in \eqref{1} with $\frac{K}{N}$ to get
\begin{align}\label{inequality}
D_t\leq&\beta^2D_{t-1}+\beta^2\frac{K^2}{N^2}\left(1+2\sum_{j=2}^{t-1}\beta^{t-j}\right)+\alpha^2\frac{\sigma^2}{2\theta}+2\alpha\beta^{t}\gamma(t-1)+2\alpha^2\sum_{k=1}^{t-1}\beta^{k}\gamma(k)\\
=&\beta^2D_{t-1}+\beta^2\frac{1+\beta-2\beta^{t-1}}{1-\beta}\frac{K^2}{N^2}+\alpha^2\frac{\sigma^2}{2\theta}+2\alpha\beta^{t}\gamma(t-1)+2\alpha^2\sum_{k=1}^{t-1}\beta^{k}\gamma(k)\nonumber.
\end{align}
We now replace $\gamma(k)$ with $\frac{\sigma^2}{2\theta}e^{-\frac{\theta}{N}k}$; this gives
\begin{align*}
\sum_{k=1}^{t-1}\beta^k\gamma(k)=\frac{\sigma^2}{2\theta}\sum_{k=1}^{t-1}\beta^ke^{-\frac{\theta}{N}k}=\frac{\sigma^2}{2\theta}\sum_{k=1}^{t-1}(\beta e^{-\frac{\theta}{N}})^k=\frac{\sigma^2}{2\theta}\frac{\beta e^{-\frac{\theta}{N}}-(\beta e^{-\frac{\theta}{N}})^{t}}{1-\beta e^{-\frac{\theta}{N}}},
\end{align*}
and hence
\begin{align*}
D_t\leq&\beta^2D_{t-1}+\beta^2\frac{1+\beta-2\beta^{t-1}}{1-\beta}\frac{K^2}{N^2}+\alpha^2\frac{\sigma^2}{2\theta}+\alpha\frac{\sigma^2}{\theta}\beta^{t}e^{-\frac{\theta}{N}(t-1)}+\alpha^2\frac{\sigma^2}{\theta}\frac{\beta e^{-\frac{\theta}{N}}-(\beta e^{-\frac{\theta}{N}})^{t}}{1-\beta e^{-\frac{\theta}{N}}}\\
=&\beta^2D_{t-1}+\beta^2\frac{1+\beta-2\beta^{t-1}}{1-\beta}\frac{K^2}{N^2}+\alpha^2\frac{\sigma^2}{2\theta}+\alpha\frac{\sigma^2}{\theta}e^{\frac{\theta}{N}}(\beta e^{-\frac{\theta}{N}})^t+\alpha^2\frac{\sigma^2}{\theta}\frac{\beta e^{-\frac{\theta}{N}}-(\beta e^{-\frac{\theta}{N}})^{t}}{1-\beta e^{-\frac{\theta}{N}}}.
\end{align*}
Simple algebraic manipulations on the last three terms above give
\begin{align*}
&\alpha^2\frac{\sigma^2}{2\theta}+\alpha\frac{\sigma^2}{\theta}e^{\frac{\theta}{N}}(\beta e^{-\frac{\theta}{N}})^t+\alpha^2\frac{\sigma^2}{\theta}\frac{\beta e^{-\frac{\theta}{N}}-(\beta e^{-\frac{\theta}{N}})^{t}}{1-\beta e^{-\frac{\theta}{N}}}\\
&\quad\quad=\alpha\frac{\sigma^2}{\theta}\frac{e^{\frac{\theta}{N}}-1}{1-\beta e^{-\frac{\theta}{N}}}(\beta e^{-\frac{\theta}{N}})^t+\alpha^2\frac{\sigma^2}{2\theta}\frac{e^{\frac{\theta}{N}}+\beta}{e^{\frac{\theta}{N}}-\beta};
\end{align*}
this implies
\begin{align*}
	D_t\leq&\beta^2D_{t-1}+\beta^2\frac{1+\beta-2\beta^{t-1}}{1-\beta}\frac{K^2}{N^2}+\alpha\frac{\sigma^2}{\theta}\frac{e^{\frac{\theta}{N}}-1}{1-\beta e^{-\frac{\theta}{N}}}(\beta e^{-\frac{\theta}{N}})^t+\alpha^2\frac{\sigma^2}{2\theta}\frac{e^{\frac{\theta}{N}}+\beta}{e^{\frac{\theta}{N}}-\beta}\\
	=&\beta^2D_{t-1}+\beta^2\frac{1+\beta}{1-\beta}\frac{K^2}{N^2}+\alpha^2\frac{\sigma^2}{2\theta}\frac{e^{\frac{\theta}{N}}+\beta}{e^{\frac{\theta}{N}}-\beta}-\frac{2}{1-\beta}\frac{K^2}{N^2}\beta^{t+1}\\
	&+\alpha\frac{\sigma^2}{\theta}\frac{e^{\frac{\theta}{N}}-1}{1-\beta e^{-\frac{\theta}{N}}}(\beta e^{-\frac{\theta}{N}})^t\\
	=&\beta^2D_{t-1}+C_1-C_2\beta^{t+1}+C_3(\beta e^{-\frac{\theta}{N}})^t,
\end{align*}
with
\begin{align}\label{C}
	C_1:=\beta^2\frac{1+\beta}{1-\beta}\frac{K^2}{N^2}+\alpha^2\frac{\sigma^2}{2\theta}\frac{e^{\frac{\theta}{N}}+\beta}{e^{\frac{\theta}{N}}-\beta},\quad C_2:=\frac{2}{1-\beta}\frac{K^2}{N^2},\quad C_3:=\alpha\frac{\sigma^2}{\theta}\frac{e^{\frac{\theta}{N}}-1}{1-\beta e^{-\frac{\theta}{N}}}.
\end{align}
Solving for $t\in\{2,...,N\}$ the following recursive inequality
\begin{align*}
	D_t\leq&\beta^2D_{t-1}+C_1-C_2\beta^{t+1}+C_3(\beta e^{-\frac{\theta}{N}})^t,
\end{align*}
with $D_1=\frac{\sigma^2}{2\theta}$, yields
\begin{align*}
D_N\leq& \beta^{2(N-1)}\frac{\sigma^2}{2\theta}+\sum_{t=2}^N\beta^{2(N-t)}(C_1-C_2\beta^{t+1}+C_3(\beta e^{-\frac{\theta}{N}})^t)\\
=&\beta^{2(N-1)}\frac{\sigma^2}{2\theta}+C_1\frac{1-\beta^{2(N-1)}}{1-\beta^2}-C_2\beta^{2N}\sum_{t=2}^N(\beta^{-1})^{t-1}+C_3\beta^{2N}\sum_{t=2}^N(\beta^{-1}e^{-\frac{\theta}{N}})^t\\
=&\beta^{2(N-1)}\frac{\sigma^2}{2\theta}+C_1\frac{1-\beta^{2(N-1)}}{1-\beta^2}-C_2\frac{\beta^{2N}-\beta^{N+1}}{\beta-1}+C_3\frac{\beta^{2N-1}e^{-\frac{\theta}{N}}-\beta^N e^{-\theta}}{\beta e^{\frac{\theta}{N}}-1}.
\end{align*}
In conclusion, recalling notation \eqref{C}, we obtain, for all $N\geq 2$, the bound
\begin{align*}
D_N\leq&\beta^{2(N-1)}\frac{\sigma^2}{2\theta}+\left(\beta^2\frac{1+\beta}{1-\beta}\frac{K^2}{N^2}+\alpha^2\frac{\sigma^2}{2\theta}\frac{e^{\frac{\theta}{N}}+\beta}{e^{\frac{\theta}{N}}-\beta}\right)\frac{1-\beta^{2(N-1)}}{1-\beta^2}\nonumber\\
&-\frac{2}{1-\beta}\frac{K^2}{N^2}\frac{\beta^{2N}-\beta^{N+1}}{\beta-1}+\alpha\frac{\sigma^2}{\theta}\frac{e^{\frac{\theta}{N}}-1}{1-\beta e^{-\frac{\theta}{N}}}\frac{\beta^{2N-1}e^{-\frac{\theta}{N}}-\beta^N e^{-\theta}}{\beta e^{\frac{\theta}{N}}-1},
\end{align*}
with $D_1=\frac{\sigma^2}{2\theta}$: this proves formula \eqref{upper bound}. The statement regarding the sharpness of the inequality is immediately obtained by observing that the use of an inequality bound throughout the proof occurs only at \eqref{inequality}, where we upper bound each $K_j$ with $\frac{K}{N}$. Recalling that $K_j:=m^{\star}_j-m^{\star}_{j-1}$ and that $m^{\star}_j=m^{\star}\left(\frac{j}{N}\right)$, we see that in the case of a linear drift $m^{\star}:[0,1]\to\mathbb{R}$ we can write
\begin{align*}
K_j=m^{\star}_j-m^{\star}_{j-1}=m^{\star}\left(\frac{j}{N}\right)-m^{\star}\left(\frac{j-1}{N}\right)=\frac{K}{N},
\end{align*}  	
obtaining the sharpness.
	
\subsection{Proof of formula \eqref{LLN}}\label{proof LLN}

We have
\begin{align*}
\mathbb{E}\left[\left(\frac{1}{N}\sum_{i=1}^NX_{\frac{i}{N}}-m^{\star}\right)^2\right]=&\mathbb{E}\left[\left(\frac{1}{N}\sum_{i=1}^N\left(X_{\frac{i}{N}}-m^{\star}\right)\right)^2\right]\\
=&\frac{1}{N^2}\sum_{i,j=1}^N\mathbb{E}\left[\left(X_{\frac{i}{N}}-m^{\star}\right)\left(X_{\frac{j}{N}}-m^{\star}\right)\right]\\
=&\frac{1}{N^2}\sum_{i,j=1}^N\mathtt{cov}\left(X_{\frac{i}{N}},X_{\frac{j}{N}}\right).
\end{align*}
We now utilize \eqref{12} to obtain
\begin{align*}
\mathbb{E}\left[\left(\frac{1}{N}\sum_{i=1}^NX_{\frac{i}{N}}-m^{\star}\right)^2\right]=&\frac{1}{N^2}\sum_{i,j=1}^N\mathtt{cov}\left(X_{\frac{i}{N}},X_{\frac{j}{N}}\right)\\
=&\frac{1}{N^2}\sum_{i,j=1}^N\frac{\sigma^2}{2\theta}e^{-\frac{\theta}{N}|i-j|}\\
=&\frac{1}{N^2}\frac{\sigma^2}{2\theta}\left(N+2\sum_{i>j}e^{-\frac{\theta}{N}(i-j)}\right)\\
=&\frac{1}{N^2}\frac{\sigma^2}{2\theta}\left(N+2\sum_{i=2}^N\sum_{j=1}^{i-1}e^{-\frac{\theta}{N}(i-j)}\right)\\
=&\frac{1}{N^2}\frac{\sigma^2}{2\theta}\left(N+2\sum_{i=2}^Ne^{-\frac{\theta}{N}i}\sum_{j=1}^{i-1}e^{\frac{\theta}{N}j}\right)\\
=&\frac{1}{N^2}\frac{\sigma^2}{2\theta}\left(N+2\sum_{i=2}^Ne^{-\frac{\theta}{N}i}\frac{e^{\frac{\theta}{N}}-e^{\frac{\theta}{N}i}}{1-e^{\frac{\theta}{N}}}\right)\\
=&\frac{1}{N^2}\frac{\sigma^2}{2\theta}\left(N+\frac{2}{1-e^{\frac{\theta}{N}}}\sum_{i=2}^N\left(e^{-\frac{\theta}{N}(i-1)}-1\right)\right)\\
=&\frac{1}{N^2}\frac{\sigma^2}{2\theta}\left(N+\frac{2}{1-e^{\frac{\theta}{N}}}\left(\frac{e^{-\frac{\theta}{N}}-e^{-\theta}}{1-e^{-\frac{\theta}{N}}}-(N-1)\right)\right)\\
=&\frac{1}{N^2}\frac{\sigma^2}{2\theta}\frac{Ne^{-\frac{\theta}{N}}-Ne^{\frac{\theta}{N}}+2-2e^{-\theta}}{(1-e^{\frac{\theta}{N}})(1-e^{-\frac{\theta}{N}})}.
\end{align*} 	
Comparing the first and last members above we get the desired formula \eqref{LLN}.
	
\subsection{Proof of formula \eqref{MLE}}\label{proof MLE}

Adding and subtracting $m^{\star}$ from $X_{\frac{j}{N}}$ in
\begin{align*}
\hat{M}=\frac{\sum\limits_{i,j=1}^Na_{ij}X_{\frac{j}{N}}}{\sum\limits_{i,j=1}^Na_{ij}}
\end{align*}
gives
\begin{align*}
\hat{M}=\frac{\sum\limits_{i,j=1}^Na_{ij}\left(X_{\frac{j}{N}}-m^{\star}\right)+\sum\limits_{i,j=1}^Na_{ij}m^{\star}}{\sum\limits_{i,j=1}^Na_{ij}}=\frac{\sum\limits_{i,j=1}^Na_{ij}\left(X_{\frac{j}{N}}-m^{\star}\right)}{\sum\limits_{i,j=1}^Na_{ij}}+m^{\star},
\end{align*}
and hence
\begin{align*}
\hat{M}-m^{\star}=\frac{\sum\limits_{i,j=1}^Na_{ij}\left(X_{\frac{j}{N}}-m^{\star}\right)}{\sum\limits_{i,j=1}^Na_{ij}}.
\end{align*}
We can therefore compute the mean square error of the (global) maximum likelihood estimator as 
\begin{align*}
&\mathbb{E}[(\hat{M}-m^{\star})^2]=\frac{\sum\limits_{i,j=1}^N\sum\limits_{k,l=1}^Na_{ij}a_{kl}\mathbb{E}\left[\left(X_{\frac{j}{N}}-m^{\star}\right)\left(X_{\frac{l}{N}}-m^{\star}\right)\right]}{\left(\sum\limits_{i,j=1}^Na_{ij}\right)^2}\\
&\quad\quad=\frac{\sum\limits_{i,j=1}^N\sum\limits_{k,l=1}^Na_{ij}a_{kl}c_{jl}}{\left(\sum\limits_{i,j=1}^Na_{ij}\right)^2}=\frac{\sum\limits_{i=1}^N\sum\limits_{k,l=1}^Na_{kl}\sum\limits_{j=1}^Na_{ij}c_{jl}}{\left(\sum\limits_{i,j=1}^Na_{ij}\right)^2}\\
&\quad\quad=\frac{\sum\limits_{i=1}^N\sum\limits_{k,l=1}^Na_{kl}\delta_{il}}{\left(\sum\limits_{i,j=1}^Na_{ij}\right)^2}=\frac{\sum\limits_{i=1}^N\sum\limits_{k=1}^Na_{ki}}{\left(\sum\limits_{i,j=1}^Na_{ij}\right)^2}\\
&\quad\quad=\frac{1}{\sum\limits_{i,j=1}^Na_{ij}}.
\end{align*}
Here, $\{\delta_{il}\}_{i,l=1,...,N}$ stands for the Kronecker symbol and we have exploited the fact that $A$ is the inverse of $C$. On the other hand, a direct verification shows that
\begin{align*}
A=\frac{2\theta}{\sigma^2(1-e^{-2\frac{\theta}{N}})}
\begin{bmatrix}
1 & -e^{-\frac{\theta}{N}} & 0 & 0 & \cdot & 0 \\
-e^{-\frac{\theta}{N}} & 1+e^{-2\frac{\theta}{N}} & -e^{-\frac{\theta}{N}} & 0 &  \cdot & 0\\
0 & -e^{-\frac{\theta}{N}} & 1+e^{-2\frac{\theta}{N}} & -e^{-\frac{\theta}{N}} &  \cdot & 0\\
\cdot & \cdot & \cdot & \cdot & \cdot & \cdot \\
\cdot & \cdot & \cdot & \cdot & \cdot & \cdot \\
\cdot & \cdot & \cdot & \cdot & \cdot & \cdot \\
0 & \cdot & 0 & -e^{-\frac{\theta}{N}} & 1+e^{-2\frac{\theta}{N}} & -e^{-\frac{\theta}{N}} \\
0 & 0 & \cdot & 0 & -e^{-\frac{\theta}{N}} & 1
\end{bmatrix},
\end{align*}
and hence
\begin{align*}
\sum\limits_{i,j=1}^Na_{ij}=\frac{2\theta}{\sigma^2(1-e^{-2\frac{\theta}{N}})}\left(N+(N-2)e^{-2\frac{\theta}{N}}-2(N-1)e^{-\frac{\theta}{N}}\right).
\end{align*}
To conclude,
\begin{align*}
\mathbb{E}[(\hat{M}-m^{\star})^2]=\frac{\sigma^2(1-e^{-2\frac{\theta}{N}})}{2\theta\left(N+(N-2)e^{-2\frac{\theta}{N}}-2(N-1)e^{-\frac{\theta}{N}}\right)}.
\end{align*}

\newpage

\bibliographystyle{abbrv}
\bibliography{mean_reverting_OU}

\end{document}